\documentclass[reqno, 11pt]{amsart}
\usepackage{fullpage}
\selectfont
\usepackage{algorithm,algpseudocode,mathtools}

\algnewcommand{\IIf}[1]{\State\algorithmicif\ #1\ \algorithmicthen}
\algnewcommand{\EndIIf}{\unskip\ \algorithmicend\ \algorithmicif}
\usepackage[utf8]{inputenc}
\usepackage[OT2,T1]{fontenc}
\usepackage[english]{babel}

\usepackage{amsmath}
\usepackage{amsfonts}
\usepackage{amssymb}
\usepackage{amsthm}

\usepackage{mathrsfs}
\usepackage{listings}
\usepackage[mathcal]{euscript}
\usepackage{bbm}
\usepackage{enumitem}
\usepackage{multirow}

\usepackage[dvipsnames]{xcolor}
\usepackage{graphicx}

\usepackage{tikz-cd}
\usetikzlibrary{matrix,positioning}
\usepackage[all]{xy}

\usepackage{colortbl}
\definecolor{mygray}{gray}{0.92}

\usepackage{array}
\newcolumntype{C}[1]{>{\centering\arraybackslash$}p{#1}<{$}}

\usepackage{breqn}
\newcounter{myequation}[equation]

\usepackage{float}
\restylefloat{table}
\usepackage{multirow}

\usepackage{hyperref}

\theoremstyle{plain}
\newtheorem{theorem}{Theorem}

\newtheorem{proposition}[theorem]{Proposition}

\theoremstyle{definition}

\theoremstyle{remark}
\newtheorem{remark}[theorem]{Remark}

\numberwithin{equation}{section}

\def\epsilon{\varepsilon}

\def\tilde{\widetilde}

\def\Magma{\textsc{Magma}}

\DeclareMathOperator{\Aut}{Aut}

\DeclareMathOperator{\End}{End}

\DeclareMathOperator{\id}{Id}

\DeclareMathOperator{\Sp}{Sp}

\DeclareMathOperator{\Th}{\textrm{Th}_2}

\DeclareMathOperator{\irr}{irr}
\DeclareMathOperator{\covgon}{cov.gon}

\newcommand{\dq}{{/\kern -3pt/}}

\def\C{\mathbb{C}}

\def\F{\mathbb{F}}

\def\P{\mathbb{P}\,}
\def\Q{\mathbb{Q}}

\usepackage{bm}

\def\Ac{\mathcal{A}}

\def\Hc{\mathcal{H}}

\def\Mc{\mathcal{M}}

\newcommand{\CC}{{\mathbb{C}}}

\newcommand{\HH}{{\mathbb{H}}}

\newcommand{\PP}{{\mathbb{P}}}
\newcommand{\QQ}{{\mathbb{Q}}}

\newcommand{\ZZ}{{\mathbb{Z}}}

\newcommand{\LL}{{\mathbb{L}}}

\newcommand{\Ssing}{S^{\operatorname{sing}}}

\hypersetup{
  pdfauthor   = {Samuel Grushevsky, Christophe Ritzenthaler},
  pdftitle    = {A complete curve of genus 5 in the moduli space of curves of genus 3 and related questions},
  pdfsubject  = {},
  pdfkeywords = {},
  backref=true, pagebackref=true, hyperindex=true, colorlinks=true,
  breaklinks=true, urlcolor=blue, linkcolor=blue, citecolor=blue,
  bookmarks=true, bookmarksopen=true
}

 \usepackage{color}

\begin{document}

\title{A complete curve of genus 5 in the moduli space of curves of genus $3$ and related questions }
\date{\today}

\begin{abstract}
Let $\Mc_3$ be the coarse moduli space of smooth complex projective curves of genus $3$. Through a general point of $\Mc_3$, we construct a smooth complete curve of genus 21 and gonality 7, and a singular complete curve whose normalization has genus 19 and gonality 6. We also construct a Hodge-generic complete trigonal curve of geometric genus 5 in $\Mc_3$. 
\end{abstract}

\author[S. Grushevsky]{Samuel Grushevsky}
\address{Department of Mathematics and Simons Center for Geometry and Physics, Stony Brook University, Stony Brook, NY 11794-3651, USA}

\email{samuel.grushevsky@stonybrook.edu}

\author[Ch. Ritzenthaler]{Christophe Ritzenthaler}
\address{Christophe Ritzenthaler,
  Univ Rennes, CNRS, IRMAR - UMR 6625, F-35000
 Rennes, France. }
  
\address{Christophe Ritzenthaler,
  Université Côte d'Azur, CNRS, LJAD UMR 7351,
  Nice,
  France
}
\email{christophe.ritzenthaler@univ-rennes1.fr}

\dedicatory{To Riccardo Salvati Manni, with admiration for his intuition on thetas}


\maketitle
We work over the field of complex numbers $\CC$. For $g \geq 2$, let $\Mc_g $ be the {\em coarse} moduli space of smooth projective curves of genus $g$ over $\CC$. It is a quasi-projective variety of dimension $3g-3$ and it is never complete.  In \cite{harris-complete}, Harris asked about the largest dimension of a \emph{projective}
subvariety of $\Mc_g$. This is still wide open as soon as $g \geq 4$. It is known that it is less than or equal to $g-2$ \cite{diaz} (and therefore is 1 for $g=3$ and at most $2$ for $g=4$). 

It is however easy to prove that there exist complete curves in $\Mc_g$ for any $g \geq 3$. Indeed, for this range of~$g$, the boundary of the projective Satake compactification $\Mc_g^{\textrm{Sat}}\subset\PP^N$ of $\Mc_g$ has codimension 2. Therefore, intersecting $\Mc_g^{\textrm{Sat}}$ with general hypersurfaces of sufficiently high degree in $\PP^N$ gives the existence of a curve passing through any finite collection of points of $\Mc_g$ and avoiding the boundary $\Mc_g^{\textrm{Sat}}\setminus\Mc_g$. Still, Harris then adds `It should be said, however, that
despite the apparent ubiquity of such curves, no one has yet written down
explicitly a complete curve in, for example, $\Mc_3 $.' 

Several constructions starting with \cite{kodaira-complete, kas-complete, ATI69,HIR69,riera} had been given for complete curves in $\Mc_g$ for $g \geq 4$; many of these are contained in the loci of curves with a non-trivial automorphism.
Kodaira's construction in  fact gives a complete family of smooth curves of genus 6 over a base curve of genus 9 \cite[Thm.~1.17]{zaal-phd}; the construction of Gonz\'alez-D\'iez and Harvey \cite{complete-harvey} yields a complete family of smooth curves of genus 4 over a base curve of genus 9 \cite[Thm.~1.19]{zaal-phd}; more recently for any $g>3$, Stewart \cite{stewart} generalizes the work of \cite{complete-harvey} and gives a formula for the genus of a base curve. 

What about $g=3$? The work of \cite{zaal} gives a beautiful geometric construction of an explicit complete curve in $\Mc_3$ by using a suitable covering construction. Indeed, he considers a degree 64 (possibly disconnected) cover of a genus 9 curve, but the construction does not allow one to easily compute the connected components of the curve constructed. If the cover were connected, it would be an \'etale cover of degree 64, and genus 513. One can then wonder whether an example with a precise and lower genus can be given, completing the list of examples in every genus $ \geq 3$.
In this context, we obtain the following examples (recall that the geometric genus of a singular curve is the genus of its normalization).
\begin{theorem}\label{thm:main}
\begin{enumerate}
    \item \label{item:1} There exists a smooth complete curve of genus $21$ and gonality $7$ through a general point of $\Mc_3$.
    \item \label{item:2} There exists a singular complete curve of geometric genus $19$, and whose normalization has gonality  $6$, through a general point of $\Mc_3$.
    \item\label{item:3} There exists a complete curve in $\Mc_3$ which is Hodge-generic, and whose normalization is a trigonal genus 5 curve.
\end{enumerate}
\end{theorem}
Here we recall that Hodge-genericity means that the Hodge group of a very general point is the full symplectic group. In our case this is equivalent to saying that the curve is not contained in any proper Shimura subvariety of $\Ac_3$. We note that most prior constructions of complete subvarieties of $\Mc_g$ use covering constructions and are not Hodge-generic, unlike our explicit geometric genus 5 example.

Our constructions utilize the explicit description of the Satake compactification $\Mc_3(2,4)^{\textrm{Sat}}$ of the level $(2,4)$ cover $\Mc_3(2,4)$ of the moduli space, which in the case of $g=3$ is the same as the Satake compactification $\Ac_3(2,4)^{\textrm{Sat}}$ of the level cover of the moduli space $\Ac_3$ of principally polarized complex abelian threefolds. 
As we review below, theta constants of the second order define an embedding of $\Ac_3(2,4)^{\textrm{Sat}}$ as a highly singular degree 16 hypersurface~$S\subset\PP^7$. 

The article \cite{irrat} studies various measures of irrationality of an arbitrary $n$-dimensional projective variety $X$, in particular the degree of irrationality $\irr(X)$, which is the minimal degree of a dominant rational map $X\dashrightarrow \PP^n$, and the covering gonality $\covgon(X)$, which is the minimal gonality of a normalization of a curve through a very general point of~$X$. By projection from a point, it is immediate to see that $\irr(S)\le 15$, but in fact the point $p_0\coloneqq (1:0:\ldots:0)$ has multiplicity $7$ on~$S$, so that projecting from $p_0$ gives $\irr(S)\le 9$. We will further improve this by an explicit computation, obtaining
\begin{proposition}\label{prop:irr}
    The degree of irrationality $\irr(S)\le 8$. 
\end{proposition}
We note that in general the degree of irrationality gives an upper bound for the covering gonality, so the above implies $\covgon(S)\le 8$. As $\Mc_3$  itself (without a level cover) is known to be rational by \cite{katsylo}, we of course have $\covgon(\Mc_3)=\irr(\Mc_3)=1$. This, however, does not shed any light on {\em complete} curves in $\Mc_3$, as opposed to those contained in its compactification.

\begin{remark}
    The article \cite{irrat} gives the values of various measures of irrationality for very general hypersurfaces in $\PP^n$, under various conditions on the degree and dimension. However, most of their results apply to hypersurfaces that are smooth or have at worst canonical singularities. The singularity of $S$ at the point $p_0$ is not canonical, and thus the lower bounds from \cite{irrat} do not seem to apply directly.
\end{remark}

We conclude this introduction with two further remarks on related questions.
\begin{remark}
    
While for the {\em existence} of complete curves, it does not matter whether one considers the coarse moduli space $\Mc_g$ or the stack $M_g$, if one is interested in determining the {\em genera} or other invariants of complete curves, working in $\Mc_g$ is completely different from working in $M_g$. Oort's work \cite{OOR74} explains this crucial distinction. Indeed, while $M_g$ has a metric of negative curvature, and thus does not contain any complete rational or elliptic curves, Oort shows that there exists an infinite collection of genera for which $\Mc_g$ contains a complete rational curve (see also \cite{debarre}).

The data of a nonconstant morphism from a smooth complete curve to the stack $M_g$ is equivalent to the data of a non-isotrivial smooth family over that curve.  Bryan and Donagi \cite{BD02} constructed such surface bundles, and in particular showed that for an infinite collection of genera of the form $2(2n-1)n^2+1$ for $n\ge2$ (so starting with $M_{25}$) the moduli stack $M_g$ contains a complete genus 2 curve. While there are some further bounds for the minimal genus of a complete curve in $M_g$ coming from surface theory (see e.g. \cite{KOT98}), no lower bounds for the genera of complete curves in $\Mc_g$ seem to be known. We refer to \cite{CAT17} for a survey, and to \cite{FLA22} for some specific considerations in genus 3.
\end{remark}

\begin{remark}
    It is also natural to ask about the minimal genus and gonality of complete curves in (the coarse space) $\Ac_g$. While we do not pursue the question of the minimal genus of a complete curve through a general point of $\Ac_g$, we note that for any $g\ge 2$ there exist complete rational curves in $\Ac_g$. Indeed, while the moduli space of hyperelliptic curves $\Hc_g\simeq\Mc_{0,2g+2}/S_{2g+2}$ is affine,
 and thus does not contain any complete curves, one can construct one-parameter families of hyperelliptic curves that only have degenerations of compact type, which is ensured by getting the branch points to only coincide in triples. One particular example of such a family is 
    \[ y^2=F_t(x), \; \textrm{with } F_t(x)=\frac{x^{2g+1}}{2g+1}-\frac{2x^{g+1}}{g+1}+x-t\,,\]
    with parameter $t\in\PP^1$. Indeed, this polynomial $F_t(x)$ is chosen precisely in such a way that its derivative is $F_t'(x)=(x^g-1)^2$.
    Thus any $t$ for which $F_t(x)$ has a multiple root has  to satisfy $$t=\frac{\zeta^{2g+1}}{2g+1}-\frac{2\zeta^{g+1}}{g+1}+\zeta= c \zeta$$ for any $g$'th root of unity $\zeta$ and $c \in \C$. For such values of $t$,
    $F_t'(x)$ then vanishes with multiplicity 2 at $\zeta$, so that $\zeta$ is then precisely a triple root of $F_t(x)$. The stable reduction of $y^2=F_t(x)$ is therefore of compact type for any finite value of $t$ (see for instance \cite[Th.5.5]{userguide}). One also verifies that the fiber at $t=\infty$ is given by the value at $v=0$ of $$Y^2 = \frac{X^{2g+1}}{2g+1} - \frac{2}{g+1}\,v^{2g}X^{g+1} \;+\; v^{4g}X -1,$$ with the change of variable $t=v^{-4g-2}, x=v^{-2} X$ and $y=v^{-2g-1} Y$. Hence it is a smooth curve and altogether the family above
    gives a complete rational curve in $\Ac_g$. This is not a trivial family since the general fiber is indecomposable whereas the ppav at $t=c \zeta$ is decomposable. 

    The above construction of a rational curve in $\Ac_3$ is not the first such in the literature: by taking a section of the Picard modular variety constructed by Shiga \cite{shiga}, not passing through any cusps, one can also obtain a rational curve (see also the detailed discussion in \cite{clvdg}). It turns out that the curve we constructed above is Hodge-generic, but we do not pursue this here.
\end{remark}

\subsection*{Acknowledgments.} The general strategy using the embedding $S$ was developed without any assistance of AI. Then, the second-named author used the Gemini large language model to draft arguments for the generators and representation of $\Sp(6,\ZZ)$ (before finding appropriate references). The first-named author had a long conversation with ChatGPT exploring possible surfaces in $\PP^7$ invariant under subgroups of the deck transformation group $G=\Sp(6,\ZZ)/\Gamma_3(2,4)$ of the level cover, and possible singularities of curves obtained as intersections of $S$ with such surfaces. This resulted in ChatGPT finding the lowest genus examples in the main theorem, the validity of all of which was verified independently by the authors, not using AI. All of the text of this paper was written by the authors, not by AI systems, and the authors of course take full responsibility for all the mathematical content.

The second-named author thanks David Lubicz and Damien Robert for references concerning theta functions and moduli spaces.

\section{Review on genus 3 theta constants}

We recall the usual notation for \emph{theta constants} with integral characteristics:
\begin{equation*}
    \theta\begin{bmatrix} \epsilon \\ \delta \end{bmatrix}(\tau) \coloneqq \sum_{m \in \ZZ^g} \exp \pi i \left[ \left(m + \frac{\epsilon}{2}\right)^\top \tau \left(m + \frac{\epsilon}{2}\right) + 2 \, \left(m + \frac{\epsilon}{2}\right)^\top \frac{\delta}{2} \right]\,.
\end{equation*}
Here $\epsilon, \delta$ are vectors of length $g$ consisting of zeroes and ones, and by parity the above theta constant vanishes identically unless $\epsilon^\top \delta \equiv 0 \pmod 2$.

Similarly the \emph{theta constants of the second order} are defined for $\sigma\in(\ZZ/2\ZZ)^g$  as
\begin{equation*} 
    \Theta[\sigma](\tau) \coloneqq \sum_{m \in \ZZ^g} \exp 2 \pi i \left[\left(m + \frac{\sigma}{2}\right)^\top \tau \left(m + \frac{\sigma}{2}\right)  \right] = \theta\begin{bmatrix} \sigma \\ 0 \end{bmatrix}(2\tau)\,.
\end{equation*}
These two types of theta constants are related by Riemann's bilinear relations
\begin{equation} \label{eq:dupli}
    \theta^2\begin{bmatrix} \epsilon \\ \delta \end{bmatrix}(\tau) = \sum_{\sigma}  (-1)^{\sigma\cdot\delta} \Theta[\sigma + \epsilon](\tau) \cdot \Theta[\sigma](\tau)\,. \tag{1}
\end{equation}
In the above formulas, $\tau\in\HH_g$ is a matrix in the Siegel upper half space, i.e.~it is a symmetric complex $g\times g$ matrix with positive definite imaginary part. The integral symplectic group $\Sp(2g,\ZZ)$ acts on $\HH_g$, and the quotient $\Ac_g\coloneqq \HH_g/\Sp(2g,\ZZ)$ is the moduli space of complex principally polarized abelian varieties (ppav) of dimension~$g$.

Theta constants with characteristics, and theta constants of the second order, are~{\em not} Siegel modular forms with respect to the full symplectic group $\Sp(2g,\ZZ)$, that is they are not sections of a line bundle on $\Ac_g$.  However, theta constants {\em are} modular forms with respect to suitable theta level subgroups, which we now recall. For an even $n>0$, the level $(n,2n)$ theta subgroup is defined as
$$
\Gamma_g(n,2n)\coloneqq \left\{ \begin{pmatrix}
    A & B \\ C & D 
\end{pmatrix} \in \Sp(2g,\ZZ),\; 
\begin{array}{c}
     A\equiv D \equiv \id_g \pmod{n}  \\
      B \equiv C  \equiv 0 \pmod{n} \\
      \mathrm{diag} A^\top B \equiv \mathrm{diag} C^\top D \equiv 0 \pmod{2n}
\end{array} \right\}.$$

Theta constants with characteristics are modular forms of weight 1/2 with respect to $\Gamma_g(4,8)$, while theta constants of the second order are modular forms, also of weight $1/2$, with respect to the bigger subgroup $\Gamma_g(2,4)$. We will not recall the full theta transformation formula or even the notion of Siegel modular forms (and refer to \cite{igusa1} for all of the above), but record that this means that there is a well-defined level two theta map
$$
 \Th:\Ac_g(2,4)\to\PP^{2^g-1}
$$
which is induced by the map sending a point $\tau\in\HH_g$ to the set of values of all theta constants of the second order $\{\Theta[\sigma](\tau)\}_{\sigma\in(\ZZ/2\ZZ)^g}$; here $\Ac_g(2,4)\coloneqq \HH_g/\Gamma_g(2,4)$ is the so-called level moduli space of ppav, which is a finite cover of~$\Ac_g$. We will denote by $\pi:\Ac_g(2,4)\to\Ac_g$ this covering map, and denote $G\coloneqq\Sp(2g,\ZZ)/\Gamma_g(2,4)$ its Galois group. We refer to \cite{manni-level} for many more foundational details on these level covers.

It is known (see \cite{igusa1}) that the analogous map defined by taking all theta constants with characteristics is an embedding of $\Ac_g(4,8)$, and in \cite{manni-emb} it is claimed that $\Th$ is also injective, for any $g\ge 3$. However, the argument there has a gap that has to do with signs of theta constants, and thus for arbitrary $g\ge 3$ the map $\Th$ is only known to be finite.\footnote{A proof of the actual non-injectivity of $\Th$ on $\Ac_g(2,4)$ for $g\ge 6$ will be published elsewhere.} 

This said, the specific situation in genus 3 is completely understood. We refer to \cite{glass} for the classical theta constant treatment in genus 3, to \cite[II]{complete-harvey} for the discussions of the extension of  theta constants to  the boundary of the Satake compactification, adapted to our viewpoint, and to \cite{univkummer} for a quick modern summary of relevant results. The map $\Th:\Ac_3(2,4)\to \PP^7$ {\em is} an embedding, and its image is a hypersurface, as we now review, following \cite{geemen}. Indeed, one starts with the following form of Riemann's quartic relation among genus 3 theta constants with characteristics:
\begin{align*}
    & \theta\begin{bmatrix} 0 & 0 & 0 \\ 0 & 0 & 0 \end{bmatrix}(\tau) \;
      \theta\begin{bmatrix} 0 & 0 & 0 \\ 1 & 0 & 0 \end{bmatrix}(\tau) \;
      \theta\begin{bmatrix} 0 & 0 & 0 \\ 0 & 1 & 0 \end{bmatrix}(\tau) \;
      \theta\begin{bmatrix} 0 & 0 & 0 \\ 1 & 1 & 0 \end{bmatrix}(\tau) \\
    - \, & \theta\begin{bmatrix} 0 & 0 & 1 \\ 0 & 0 & 0 \end{bmatrix}(\tau) \;
      \theta\begin{bmatrix} 0 & 0 & 1 \\ 1 & 0 & 0 \end{bmatrix}(\tau) \;
      \theta\begin{bmatrix} 0 & 0 & 1 \\ 0 & 1 & 0 \end{bmatrix}(\tau) \;
      \theta\begin{bmatrix} 0 & 0 & 1 \\ 1 & 1 & 0 \end{bmatrix}(\tau)  \\
    - \, & \theta\begin{bmatrix} 0 & 0 & 0 \\ 0 & 0 & 1 \end{bmatrix}(\tau) \;
      \theta\begin{bmatrix} 0 & 0 & 0 \\ 1 & 0 & 1 \end{bmatrix}(\tau) \;
      \theta\begin{bmatrix} 0 & 0 & 0 \\ 0 & 1 & 1 \end{bmatrix}(\tau) \;
      \theta\begin{bmatrix} 0 & 0 & 0 \\ 1 & 1 & 1 \end{bmatrix}(\tau) = 0.
\end{align*}
We want to express this as an identity among theta constants of the second order. However, this would involve square roots arising from Riemann's bilinear relations \eqref{eq:dupli}, and thus to get a polynomial equation in $\Theta[\sigma](\tau)$ we need to take Galois conjugates. This is to say, if we write the above relation as $r_1-r_2-r_3=0$, we need to multiply the four such equations, with all choices of signs of $r_2,r_3$. This yields the equation
$$
 F\coloneqq r_1^4+r_2^4+r_3^4-2r_1^2r_2^2-2r_1^2r_3^2-2r_2^2 r_3^2=0\,.
$$

For ease of notation, from now on we will write $T_0,\dots,T_7$ for the coordinates on $\PP^7$ corresponding to the theta constants of the second order: $T_i$ corresponds to $\Theta[\sigma](\tau)$ where $\sigma$ is the binary expansion of $i$, written in reverse, so that for example $T_1=\Theta[100](\tau)$ and $T_3=\Theta[110](\tau)$.
As discussed in {\it loc.~cit.}, $S\coloneqq \{F=0\}\subset\PP^7$ is an irreducible hypersurface of degree 16, which thus coincides with the closure of $\Th(\Ac_3(2,4))$, i.e. with $\Th(\Ac_3(2,4)^{\textrm{Sat}})$. Moreover, the complement is known: $S\setminus \Th(\Ac_3(2,4))$ is the union of 126 copies of $\PP^3$, each given by the vanishing of a certain configuration of 16 (even) theta constants with characteristics.

Furthermore, every principally polarized abelian threefold is the Jacobian of a smooth curve, or a product of Jacobians (see \cite{oort-ueno} for a fully general treatment of this). Thus the complement \[\Ssing\coloneqq S\setminus \Th(\textrm{Jac}(\Mc_3(2,4))),\] where $\textrm{Jac}$ is the Jacobian map,  is the closure of the locus of decomposable ppav, i.e. it is equal to the $G$-orbit of the closure of $\Th\left(\Ac_1(2,4)\times\Ac_2(2,4)\right)$, which is a union of 336 copies of $\PP^1\times\PP^3$.

In particular, $\Ssing$ has codimension two within $S$, and moreover by \cite[Lem.~3.2]{geemen} (see also \cite{univkummer} for a more modern exposition, and more details) it is known that $S\setminus\Ssing$ is smooth, which partially justifies the notation~$\Ssing$ (though we do not claim that all points of $\Ssing$ are singular on~$S$). 

The locus $\Ssing$ is the locus where at least 2 (even) theta constants with characteristics vanish simultaneously, since the Jacobian of a smooth non-hyperelliptic genus 3 curve has no vanishing (even) theta constant, while the Jacobians of hyperelliptic genus 3 curves have exactly one vanishing even theta constant. Using \eqref{eq:dupli} to write the squares of theta constants with characteristics as quadrics $R_{\epsilon,\delta}$ in the $T_i$ gives equations for components of $\Ssing$ as components of intersections of two quadrics in $\PP^7$ with $S$.

\subsection{Idea of the constructions}
Given this geometry, in Sec.~\ref{sec:explicit} we will construct explicit complete curves in $\Mc_3$. For each of them, we will first construct a complete curve $C\subset S$, as the intersection $C\coloneqq\LL\cap S$, where in the simplest case $\LL\simeq \PP^2$ is a linear subspace of $\PP^7$, but in general is an explicit complex surface in $\PP^7$ not contained in $S$. To check that $C\subset \Th(\Mc_3(2,4))$ (and then by abuse of notation we will also denote by $C$ its preimage in $\Mc_3(2,4)$), we will verify\footnote{in most cases, algorithmically, using \Magma. The corresponding programs are provided in the ancillary file available on arXiv. Specific tricks used during the computations (e.g.~working over a finite field is sometimes sufficient and quicker) are explained in the program file.}  that at every point of $C$ at most 1 even theta constant vanishes, by intersecting $C$ with two distinct $R_{\epsilon,\delta}=0$. We will always denote by $\tilde C$ the normalization of~$C$. Note that the degree of $C$ in $\PP^7$ is equal to $16$ (the degree of $S$) times the degree of $\LL$, and in particular the lower the degree of $\LL$, the lower the degree and thus the genus of $C$. Note also that forcing $C$ to be singular, for fixed degree, will decrease its geometric genus, which we will exploit. For the genus 5 example, we will also arrange the intersection $\LL\cap S$ to be reducible, and only take one component of it. 

However, the main gain, and the reason to not just use planes $\LL\simeq\PP^2$, is to try to arrange that $\LL$ (and thus $C$) is invariant under a nontrivial subgroup $H$ of the deck group $G=\Sp(6,\ZZ)/\Gamma_3(2,4)$. We will then consider the curve $D\coloneqq \pi(C)\subset\Mc_3$, where $\pi:\Mc_3(2,4)\to\Mc_3$ is the level cover. If all $G$-images of $C$ were disjoint, of course $D$ would be isomorphic to~$C$. However, if $C$ is $H$-invariant as a set, the degree of the cover $C\to D$ is at least the order of $H$, which leads to $D$ having a drastically lower genus, and vindicates suitable higher degree choices for $\LL$. Forcing $D$ to be singular will give still more flexibility in the construction.

\subsection{The explicit deck group}
For explicit computations, we  review the action of the generators of the deck group $G$.
By \cite{runge-siegel,oura}, we have the following generators and their action on $\P^7$:
\begin{itemize}
    \item The inversion $J$ corresponds to the modular transformation $\tau \mapsto -\tau^{-1}$. The representation of $J$ on the theta constants of the second order is the Hadamard matrix 
    \begin{equation}
 \begin{split}
 J &:\ (T_\sigma)_\sigma\longmapsto
      \left(\frac{1+i}{2}\right)^3 \cdot  \left(\sum_{\mu\in\F_2^3}(-1)^{\sigma^\top \mu} \cdot T_\mu\right)_\sigma.\\
 \end{split}
\end{equation}
\item  The diagonal translations $D_j$, for $j \in \{1, 2, 3\}$, correspond to the translation $\tau \mapsto \tau + B_j$, where $B_j$ is a $3 \times 3$ symmetric matrix with a $1$ at position $(j,j)$ and $0$ elsewhere. The matrix $D_j$ acts as
\begin{equation}
 \begin{split}
 D_j&:\ T_\sigma\longmapsto i^{\sigma_j} \cdot T_\sigma.\\
 \end{split}
\end{equation}
\item  The crossed translations $D_{jk}$, for $1 \le j < k \le 3$,  correspond to the translation $\tau \mapsto \tau + B_{jk}$, where $B_{jk}$ is a $3 \times 3$ symmetric matrix with $1$'s at positions $(j,k)$ and $(k,j)$, and $0$ elsewhere. Thus, $D_{jk}$ acts as 
    \begin{equation}
 \begin{split}
 D_{jk}:\ T_\sigma\longmapsto(-1)^{\sigma_j\sigma_k} \cdot T_\sigma.\\
 \end{split}
\end{equation}
\end{itemize} 
Altogether, these matrices generate a linear group of order $4 \cdot \#G= 4 \cdot 92\,897\,280$ where the center is generated by the scalar matrix $i \cdot \textrm{Id}_8$ (acting trivially on the projective space).


\section{Explicit curve constructions} \label{sec:explicit}
\subsection{A smooth complete curve of genus $21$ through a general point of $\Mc_3$}
We start by constructing a smooth curve through a general point $p\in\Mc_3$, by intersecting $S$ with a suitable complex surface in $\PP^{7}$.

As discussed above, $S$ is smooth outside the locus $\Ssing$, which is of codimension 2. Thus by Bertini's theorem, for five general hyperplanes $H_1,\dots,H_5\subset\PP^7$ the intersection $C\coloneqq H_1\cap\dots \cap H_5\cap S$ is smooth and disjoint from $\Ssing$. But then $C\subset\PP^2\cong H_1\cap\dots\cap H_5$ is a plane curve of degree 16 (as it is the intersection with the degree 16 hypersurface $S$). Thus the genus of this plane curve is $g(C)=\tfrac{(16-1)\cdot (16-2)}{2}=105$. 

However, a much better general construction is possible by intersecting $S$ with a higher degree surface that is invariant under a subgroup $H \simeq(\ZZ/2\ZZ)^3\subset G$. This explicit construction was found by ChatGPT when searching for surfaces obtained as complete intersections, that are invariant under subgroups of~$G$. 

Indeed, given a general point $p\in S$, with coordinates $(p_0:\dots:p_7)$, we first consider the three-dimensional subspace $\PP^3\simeq \PP_p\subset\PP^7$ given by the points of the form
\[x_0p_0:x_0p_1:x_2p_2:x_3p_3:x_1p_4:x_2p_5:x_1p_6:x_3p_7\]
where $(x_0:x_1:x_2:x_3)$ are homogeneous coordinates on $\PP^3$. The reason for this awkward choice is that one can check with \Magma\ that there is a subgroup $H\subset G$ that is isomorphic to $(\ZZ/2\ZZ)^3$, which acts by changing the signs of each $x_i$, modulo the overall global sign change, that preserves $\PP_p$. The restriction of the defining degree 16 equation $F$ of the hypersurface $S$ to $\PP_p$ is a degree 16 polynomial in $x_i$, which is invariant under the change of sign of each $x_i$, and is thus a degree 8 polynomial in the squares $x_i^2$. Let now $Q=\sum a_i x_i^2$ be any linear function in $x_i^2$, with all $a_i\ne 0$, but such that $\sum a_i=0$, so that $Q(1:1:1:1)=0$. Then $Q$ defines a smooth quadric surface $\QQ_p\subset \PP_p$, and we finally let $C\coloneqq \QQ_p\cap S$. 

Since $\QQ_p$ is invariant under $H$, so is $C$, and the cover $\pi|_C:C\to D$ factors through the degree 8 cover $C\to C/H$. By thinking of $C\subset\PP_p$, we see that the curve $C/H$ is the image under the map that sends each $x_i$ to $x_i^2$, which is a degree 8 cover $\PP^3\to\PP^3$. Thus $C/H$ is a plane curve of degree 8.  We verify with \Magma\ that for one particular choice of~$p\in S$ this curve is a smooth plane octic, and thus has genus 21.
Moreover, we verify that it has no further stabilizer in $G$, and that all $G/H$ images of $C$ are disjoint, so that $D=C/H=\pi(C)\subset\Mc_3$ is a smooth genus 21 curve, and, as a smooth plane octic, it has gonality 7 \cite[p.6]{Kato1991}.  Since the plane octics $C/H$ form an algebraic family parametrized by the irreducible variety of pairs $(p,Q)$, with $p \in S$ and $Q$ as above, and the properties above are open conditions in such a family, this completes the proof of item~\eqref{item:1} in the main theorem.

\subsection{A complete singular curve of geometric genus 19 through a general point of $\Mc_3$}
While the above constructs a smooth complete curve of genus 21 through a general point of $\Mc_3$, it is now also clear how to reduce the geometric genus if the curve is no longer required to be smooth. Implementing this proves item \eqref{item:2} in Theorem \ref{thm:main}.

Indeed, as before we choose a general point $p\in S$, and consider the three-dimensional linear space $\PP_p$ containing it, and the intersection $\PP_p\cap S$. As discussed before, this intersection is given by a degree 8 equation in the squares of the coordinates, that is in $x_i^2$. We now choose a plane $\LL\simeq\PP^2$ in the quotient $\PP^3$ with coordinates $y_i=x_i^2$, bitangent to the degree 8 surface defined by the equation $F$ above in terms of the~$y$ variables. We also impose that $\LL$ passes through the point $(1:1:1:1)$, which is possible since the bitangent planes form a curve in the dual space. The pullback of this bitangent plane is then a quadric in terms of the~$x$ variables, passing through $p$, and  one checks with \Magma\ that $C/H$ has precisely two ordinary nodes corresponding to the points of bitangency. Thus $C/H$ is an irreducible plane octic with two nodes, and its normalization has genus $19=21-2$, as claimed.

As for the gonality, having a line passing through a node implies that the gonality is at most 6. Using \cite[Th.~2.1]{Kato1991}, one sees that the gonality is exactly $6$.

\subsection{A complete curve in $\Mc_3$ whose normalization is a trigonal genus 5 curve}\label{sec:g5}

We now construct a trigonal genus 5 curve in $\Mc_3$, proving item \eqref{item:3} in Theorem \ref{thm:main}. For this, we will choose a specific plane
$\LL\subset\PP^7$ invariant under a different subgroup $H\simeq\ZZ/9\ZZ$ of
the deck group~$G$. The intersection of $\LL$ with $S$ will consist of a double line (which we will discard) and a
residual curve $C$ of degree~14. The reduction in genus
comes from taking the quotient by~$H$: we will show that $C/H$ is birational to a plane curve of degree 8 whose normalization is a
trigonal curve of genus 5.  We will also check that $C$ avoids $\Ssing$ and that $H$ is its full stabilizer within $G$. This explicit construction was also found by ChatGPT given the idea for what geometric properties to look for.

For this, we think of the space $\CC^8$ of unprojectivized theta coordinates $T_i$ as $(\CC^2)^{\otimes3}$, i.e.~the coordinate $T_i$ with $i = \sigma_1 + 2 \sigma_2 + 4\sigma_3$ corresponds to $e_{\sigma_1} \otimes e_{\sigma_2} \otimes e_{\sigma_3}$. Denote then $\rho\coloneqq\exp(2\pi i/36)$,
$\eta\coloneqq\rho^4$, and $\omega\coloneqq\eta^3$. All computations will take place in $\Q(\rho)$. Consider
\[
 A_3\coloneqq-\frac1{1-i}
 \begin{pmatrix}1&i\\1&-i\end{pmatrix}.
\]
The vectors
 $a\coloneqq(1,i+(1+i)\omega)$ and $b\coloneqq(1,i+(1+i)\omega^2)$
are eigenvectors of $A_3$ with eigenvalues $\omega$ and $\omega^2$,
respectively. Define
\[
 g(v_0\otimes v_1\otimes v_2)
 \coloneqq A_3v_2\otimes v_0\otimes v_1.
\]
One can check with \Magma\ that $g$ belongs to the deck group~$G$.
Since $g^3=A_3^{\otimes3}$, it has order~9; we take
$H\coloneqq\langle g\rangle$.
To specify the plane $\LL$, we set
\[
 \begin{array}{lll}
 U_0\coloneqq b\otimes a\otimes a,&U_1\coloneqq a\otimes b\otimes a,&U_2\coloneqq a\otimes a\otimes b,\\
 V_0\coloneqq a\otimes b\otimes b,&V_1\coloneqq b\otimes a\otimes b,&V_2\coloneqq b\otimes b\otimes a,
 \end{array}
 \qquad
 \begin{aligned}
 e_1&\coloneqq U_0+\eta^2U_1+\eta^4U_2,\\
 e_2&\coloneqq V_0+\eta^4V_1+\eta^8V_2,\\
 e_8&\coloneqq V_0+\eta^7V_1+\eta^5V_2.
 \end{aligned}
\]
Then $g\circ e_j=\eta^je_j$, and we choose the invariant plane $ \LL$ parametrized by  $(x,y,z) \mapsto x e_1+y e_2+z e_8$. In these coordinates $g$ acts by
$[x:y:z]\mapsto[\eta x:\eta^2y:\eta^8z]$. Restricting the defining equation
$F$ of $S$ to this plane gives
\begin{equation}\label{eq:g5-plane-section}
 F|_{\LL}=c\cdot x^2\cdot G_{14}(x,y,z),
\end{equation}
where $c$ is a non-zero complex constant, and $G_{14}$ is a homogeneous polynomial of degree~14. Thus, unlike the preceding
constructions, we will be able to discard the double line $\ell\coloneqq\{x=0\}$ and work with the residual irreducible curve $C\coloneqq\{G_{14}=0\}\subset\LL$.

To ensure completeness of $C$, we  check in \Magma\ that $C$ does not intersect the decomposable locus, i.e. that
\[
 C\cap\{R_{\epsilon,\delta}=R_{\epsilon',\delta'}=0\}
 =\varnothing
\]
for every pair of distinct even characteristics. 

The curve $C$ is $H$-invariant, since $H$ preserves $\LL$, $S$ and the line $\ell$. Conversely,
as $C$ is an irreducible plane curve of degree $14$, it spans $\LL$, so any element of $G$
stabilizing $C$ also stabilizes $\LL$. We check with \Magma\  that $\mathrm{Stab}_G(C)=\mathrm{Stab}_G(\LL)=H$.

The generator of $H$ acts on $\LL$ by $[x:y:z]\mapsto[\eta x:\eta^2y:\eta^8z]$, and thus the functions
\[
u :=  \frac{y^3}{z^3}, \qquad t :=  \frac{y^2z}{x^3}
\]
are $H$-invariant. We check with \Magma\ that the rational map
$C\to\PP^2$, $[x:y:z]\mapsto[x^3y^3 : y^2z^4 : x^3z^3]$, has degree $9=\# H$ onto its
image $C'$. Since this rational map factors through $C/H$, the induced rational map $C/H\dashrightarrow C'$ is birational.
Moreover, since $\mathrm{Stab}_G(C)=H$, the map $C/H\to D\coloneqq \pi(C)\subset \Mc_3$ is birational as
well. The curve $C'$ is an irreducible plane curve of degree $8$, hence of arithmetic genus $21$.
Its singular points are an ordinary triple point, an ordinary node, an ordinary quadruple point,
and one non-ordinary singular point with $\delta$-invariant $6$, so that the total
$\delta$-invariant is $3+1+6+6=16$. Thus, its normalization $\tilde C'$ has genus $21-16=5$.

Finally, $\tilde C'$ is trigonal. Using \Magma, the canonical ideal of $\tilde C'\subset\PP^4$ contains
exactly three independent quadrics, which rules out the hyperelliptic case. It is not generated
by them: the three quadrics cut out a surface (a rational normal scroll) rather than the
canonical curve. By Petri's theorem, a non-hyperelliptic curve of genus $5$ whose canonical
ideal is not generated by quadrics is trigonal.

We will now prove that the constructed curve is Hodge-generic in $\Ac_3$.
By \cite[Thm.~1]{ribet}, it suffices to show that the  generic geometric
endomorphism algebra of the family of abelian threefolds over $C$ is $\QQ$.
Endomorphisms of the generic fiber specialize injectively to every fiber and
remain Hodge classes along $C$. It is then enough to find a point $p\in C$ such that all elements of $\End^0(A_p)$ which remain a Hodge class to first order
along $C$ are scalar. This is weaker than exhibiting a point $p\in C$ with $\End^0(A_p)=\Q$, which would be hard to certify rigorously; instead we work at a CM point, where $\End^0(A_p)$ is known explicitly and the first-order condition reduces to linear algebra.

We take $p=[e_2]\in C$, which is fixed by $g \in G$ above and is a smooth point of both $C$
and $S$ (checked with \Magma). Since $p \in C \subset \textrm{Jac}(\Mc_3(2,4))$, we have
$A_p\simeq \textrm{Jac}(C_p)$ for a smooth curve $C_p$ of genus~$3$. Choosing a preimage
$\tau\in\HH_3$ of $p$ and composing any lift of our $g \in H$ to
$\Sp(6,\ZZ)$ with a suitable element of $\Gamma_3(2,4)$, we get a lift of
$g$ fixing $\tau$, that is, an element $\gamma=\pm\alpha$ with
$\alpha\in\Aut(C_p)$ by Torelli theorem. As $-1\in\Gamma_3(2,4)$, $\alpha$ itself maps to
$g$, so its order is a multiple of $9$; in particular $C_p$ admits an
automorphism of order $9$. By the classification of automorphism groups of
hyperelliptic and non-hyperelliptic curves of genus~$3$ (see for instance
\cite{LR11,LRRS}), this implies
\[
  C_p\simeq C_9\coloneqq\{ y^3=x^4-x\},\qquad \Aut(C_9)\simeq\ZZ/9\ZZ,
\]
so that $\alpha$ has order exactly $9$ and $\alpha(x,y)=(\xi^3x,\xi y)$ for some
primitive ninth root of unity~$\xi$. The automorphism $\alpha$ acts on the basis
$dx/y,\ dx/y^2,\ x\,dx/y^2$ of $H^0(C_9,\Omega^1)$ by $\xi^2,\xi,\xi^4$. The CM type
$\{1,2,4\}\subset(\mathbb{Z}/9\mathbb{Z})^\times$ is primitive: $\QQ(\xi^3)$ is the only
proper CM subfield of $\QQ(\xi)$ and $\{1,2,4\}$ is not equal to any of the cosets
$\{1,4,7\}$ and $\{2,5,8\}$ of $\mathrm{Gal}(\QQ(\xi)/\QQ(\xi^3))$. Therefore $A_p$ is
simple and $\End^0(A_p)=\QQ[\alpha]\simeq\QQ(\xi)=\QQ(\eta)$. In particular $\alpha$ acts
on $H_1(A_p,\mathbb{C})$ with the six primitive ninth roots of unity as simple
eigenvalues; we denote by $E_k$ the eigenline for $\eta^k$,
$k\in(\mathbb{Z}/9\mathbb{Z})^\times$.

We compute in \Magma\ that, in suitable bases, $dg_p|_{T_pC}=\eta^8$ and
$dg_p|_{T_pS}=\mathrm{diag}(\eta,\eta^2,\eta^3,\eta^5,\eta^6,\eta^8)$.
Since  the stabilizer of $\tau$ in $\Gamma_3(2,4)$ is $\{\pm1\}$,  it
acts trivially on $\mathbb{H}_3$ and the projection $\mathbb{H}_3\to\Ac_3(2,4)$ is a
local isomorphism at $\tau$. As moreover $\mathrm{Th}_2$ is an embedding,
$d(\mathrm{Th}_2)_\tau\colon T_\tau\mathbb{H}_3\to T_pS$ is an isomorphism. By the theta
transformation formula, $\mathrm{Th}_2\circ\gamma=g\circ\mathrm{Th}_2$ as maps
$\mathbb{H}_3\to\mathbb{P}^7$, and differentiating at the fixed point $\tau$ gives
\begin{equation} \label{eq:equi}    
  d(\mathrm{Th}_2)_\tau\circ d\gamma_\tau=dg_p\circ d(\mathrm{Th}_2)_\tau .
\end{equation}
On the other hand, $d\gamma_\tau$ acts on $T_\tau\mathbb{H}_3$, the space of complex
symmetric $3\times 3$ matrices, by $X\mapsto MXM^T$, where $M$ is the matrix of the
analytic representation of $\gamma=\pm\alpha$ (the sign plays no role). Identifying
$T_\tau\mathbb{H}_3$ with $\mathrm{Sym}^2W$, where $W$ is one of the two Hodge summands
of $H_1(A_p,\mathbb{C})$, the map $d\gamma_\tau$ corresponds to the natural action of
$\alpha$ on $\mathrm{Sym}^2W$. Hence by \eqref{eq:equi}, the eigenvalues of $dg_p$ on $T_pS$ are the
products $\lambda_i\lambda_j$, $i\le j$, of the eigenvalues
$\lambda_1,\lambda_2,\lambda_3$ of $\alpha$ on $W$. We denote by $W'$ the other Hodge
summand.

The relation between $\xi$ and $\eta$ is not known a priori: write $\xi=\eta^m$ with
$m\in(\mathbb{Z}/9\mathbb{Z})^\times$. 
By the computation on $H^0(C_9,\Omega^1)$, the eigenvalues of $\alpha$ on the two Hodge
summands are $\xi,\xi^2,\xi^4$ and $\xi^{-1},\xi^{-2},\xi^{-4}$. Hence, since $\xi$ and
$\xi^{-1}$ are both powers of $\eta$, the eigenvalues of $\alpha$ on $W$ are $\eta^k$ with
$k\in m\{1,2,4\}$ for some $m\in(\mathbb{Z}/9\mathbb{Z})^\times$, not known a priori.
Their products $\lambda_i\lambda_j$ are then $\eta^k$ with $k\in m\{2,3,4,5,6,8\}$, and
this set equals $\{1,2,3,5,6,8\}$ only for $m=7$. Therefore the eigenvalues of $\alpha$ on
$W$ are $\eta,\eta^5,\eta^7$, so that $W=E_1\oplus E_5\oplus E_7$ and
$W'=E_8\oplus E_4\oplus E_2$.

Since $\alpha$ preserves the polarization,  $E_k$ and  $E_l$ are orthogonal unless
$k+l\equiv 0\pmod 9$. Hence, under the polarization, we can identify $W'$ with the dual $W^\vee$ of $W$. If
$\varepsilon_k\in E_k$ ($k=1,5,7$) is an eigenbasis of $W$, the dual basis
$(\varepsilon_1^*,\varepsilon_5^*,\varepsilon_7^*)$ of $W'$ satisfies
$\varepsilon_k^*\in E_{-k}$ since $\alpha$ acts by the contragradiant representation on $W^\vee$. As we have seen in the previous paragraph, in the basis $(\varepsilon_1,\varepsilon_5,\varepsilon_7)$,
$\alpha$ acts on $W$ by $\mathrm{diag}(\eta,\eta^5,\eta^7)$, so that $d\gamma_\tau$
multiplies the entry $X_{ij}$ by $\eta^{a_i+a_j}$, with $(a_1,a_2,a_3)=(1,5,7)$. Hence
$T_pC$, being the $\eta^8$-eigenline, is spanned by
$\theta=\varepsilon_1\otimes\varepsilon_7+\varepsilon_7\otimes\varepsilon_1\in\mathrm{Sym}^2W$, the symmetric matrix whose only
non-zero entries are $\theta_{13}=\theta_{31}=1$. Viewing $\theta$ as a symmetric map
$W'\to W$ (the differential of the period map along $C$), we have
$\theta(\varepsilon_7^*)=\varepsilon_1$.

An element $r\in\End^0(A_p)$ remains a Hodge class to first order in the direction $\theta$ only if $r\circ\theta=\theta\circ r$. Writing $r=f(\alpha)$ with
$f\in\QQ[X]$, $r$ acts on $E_k$ by $f(\eta^k)$; evaluating $r\circ\theta=\theta\circ r$
on $\varepsilon_7^*$ gives $f(\eta)\,\varepsilon_1=f(\eta^2)\,\varepsilon_1$, that is
$f(\eta)=f(\eta^2)$. Thus $f(\eta)\in\QQ(\eta)$ is fixed by $\eta\mapsto\eta^2$, which
generates $\mathrm{Gal}(\QQ(\eta)/\QQ)\simeq(\mathbb{Z}/9\mathbb{Z})^\times$. This shows
that $f(\eta)\in\QQ$ and $r$ is a scalar. The generic geometric endomorphism algebra is
therefore $\QQ$, and the curve $D$ is Hodge-generic.

\subsection{A degree 8 rational map $S\dashrightarrow \PP^6$}
While $S\subset\PP^7$ has degree 16, and thus projecting $S$ to $\PP^6$ from a general point of $S$ gives a degree 15 rational map, we describe a more special map, proving Proposition \ref{prop:irr}.

First, we observe that setting $T_0=1$ gives the explicit polynomial expression
\[
  F(1,T_1,\ldots,T_7)=-512\,T_1T_2T_3T_4T_5T_6T_7+\text{h.o.t}
\]
where the higher order terms have degree at least 8 in the variables $T_1,\dots,T_7$. This implies that the point $p_0=(1:0:\dots:0)\in S$ has multiplicity precisely 7, and thus projecting $S$ from $p_0$ gives a degree 9 rational map to $\PP^6$. To improve this bound, we further reparameterize this rational map (and this computation was suggested by ChatGPT in a conversation where we tried to find maps with higher tangency at a point). More specifically, we keep $T_0=1$, and set
\[
 (T_0,T_1,T_2,T_3,T_4,T_5,T_6,T_7) =(1,tu_1,tu_2,u_3,tu_4,u_5,u_6,t),
\]
which, varying $t\in\CC^*$ and $u_i\in\CC$ gives an isomorphism to the chart $\lbrace T_0T_7\ne 0\rbrace\simeq \CC^*\times\CC^6\subset\PP^7$. We then compute 
\[  
f(t,u_1,\dots,u_6)\coloneqq F(1,tu_1,tu_2,u_3,tu_4,u_5,u_6,t)=t^4\left(A(u)+t^4B(u)+t^8C(u)\right)
\]
for suitable polynomials $A,B,C\in\QQ[u_1,\dots,u_6]$, which by a direct computation one can further check are non-zero. Since we are working in the neighborhood where $t\ne 0$, we discard the non-zero $t^4$ factor. But then, in this neighborhood, $\lbrace f=0\rbrace$ is a hypersurface in $\CC^7$ that has degree 8 in the variable~$t$, which thus has degree of irrationality at most 8.

\bibliographystyle{alpha}
\bibliography{synthbib}
\end{document}